\documentclass{amsproc}
\usepackage{graphicx}
\usepackage{amscd}
\usepackage{amsmath}
\usepackage{amsfonts}
\usepackage{amssymb}
\theoremstyle{plain}

\newtheorem{corollary}{Corollary}

\newtheorem{definition}{Definition}

\newtheorem{proposition}{Proposition}
\newtheorem{remark}{Remark}

\newtheorem{summary}{Summary}
\newtheorem{theorem}{Theorem}
\numberwithin{equation}{section}

\begin{document}
\title[Extension of Operators]{On Extension of Operators in Classes of Finitely Equivalent Banach Spaces.
Part 1}
\author{Eugene Tokarev}
\address{B.E. Ukrecolan, 33-81 Iskrinskaya str., 61005, Kharkiv-5,Ukraine}
\email{tokarev@univer.kharkov.ua}
\dedicatory{Dedicated to the memory of S. Banach.}
\begin{abstract}In the paper is considered two problems on extension of operators whose range
space for the first problem (or domain space for the second one) belongs to
the fixed class of finite equivalence, which is generated by a given Banach
space $X$. Both problems are considered in two variants: isometric and
isomorphic ones. Received a full solution of the first problem and the
solution, close to final one for the second problem.
\end{abstract}
\maketitle

\section{Introduction}

Results presented in this paper, are similar to classical results on extending
of operators (=linear homeomorphisms) that were obtained firstly in classical
articles of I. Nachbin [1], D. Goodner [2], J. Kelley [3] and were generalized
by J. Lindenstrauss [4], [5]. In all these papers (and many others) it was
investigated a question of description of a class $\mathcal{K}$ of Banach
spaces $X$ that have the so called \textit{extension property} (more precise,
the \textit{extension 'from' property}):

\textit{For any Banach space }$Y$, \textit{which contains }$X\mathit{\in
}\mathcal{K}$ \textit{as a subspace, every operator }$u:X\rightarrow
Z$\textit{, where }$Z$\textit{ is an arbitrary Banach space, can be extended
to an operator }$\widetilde{u}:Y\rightarrow Z$\textit{ (i.e., such that the
restriction }$\widetilde{u}\mid_{X}=u$\textit{) of the same norm.}

Besides, also have been studied the \textit{extension 'to' property}. A Banach
space $X$ has a such property provided:

\textit{For any pair of Banach spaces }$Y$\textit{, }$Z$\textit{ where }%
$Y$\textit{ is a subspace of }$Z$\textit{ (in symbol: }$Y\hookrightarrow
Z$\textit{) every operator }$u:Y\rightarrow X$\textit{ can be extended to an
operator }$\widetilde{u}:Z\rightarrow X$\textit{ of the same norm.}

These properties may be reformulated in such a way that there will be
evident\ directions of generalization. Let us introduce into practice the
following designation.

Let $\mathcal{K}$, $\mathcal{W}$ and $\mathcal{V}$ be some classes of Banach
spaces; $\mathsf{R}$ be a class of linear operators.

The symbol
\[
\left[  \mathsf{R}:\left(  \mathcal{K}\subset\mathcal{W}\right)
\rightarrow\mathcal{V}\right]
\]
means:

\textit{For every Banach space }$X\in\mathcal{K}$ \textit{and every Banach
space }$Y\in\mathcal{W}$ \textit{that contains }$X$\textit{ as a subspace
every operator }$u\in\mathsf{R}$, \textit{which maps }$X$\textit{ to }%
$Z$\textit{, where }$Z\in\mathcal{V}$, \textit{may be} \textit{extended to an
operator} $\widetilde{u}:Y\rightarrow Z$\textit{ of the same norm }$\left\|
\widetilde{u}\right\|  =\left\|  u\right\|  $.

Let $\mathcal{B}$ denotes the class of all Banach spaces. Then the Hahn-Banach
theorem looks like
\[
\left[  \mathsf{F}_{1}:\left(  \mathcal{B}\subset\mathcal{B}\right)
\rightarrow\mathcal{B}\right]  ,
\]
where $\mathsf{F}_{1}$ denotes a class of all operators of rank 1 (i.e., the
class of all linear functionals); the Goodner-Kelley-Nachbin result may be
written as
\[
\left[  \mathsf{L}:\left(  C(S)\subset\mathcal{B}\right)  \rightarrow
\mathcal{B}\right]  ,
\]
where $C(S)$ denotes either a space of continuous functions on extremely
disconnected compact or a class of all such spaces and $\mathsf{L}$, as usual,
denotes the class of all (linear, bounded) operators. In this notations the
problem on extension of operators may be formulated in a following way:

\textit{Let} $\mathcal{K}$, $\mathcal{W}$ \textit{are fixed. Describe a class
}$\mathcal{X}$ \textit{of Banach spaces such that either} $\left[
\mathsf{L}:\left(  \mathcal{K}\subset\mathcal{W}\right)  \rightarrow
\mathcal{X}\right]  $ \textit{or} $\left[  \mathsf{L}:\left(  \mathcal{X}%
\subset\mathcal{W}\right)  \rightarrow\mathcal{K}\right]  $.

In these formulations the class $\mathsf{L}$ of \textit{all} operators may be
substituted with the class either of all compact or weakly compact or
Radon-Nicod\'{y}m operators and so on.

In [6] it was studied the extension 'from' property in the case when
comprehensive spaces $Y$ were finitely equivalent to $X$. In a given notation
it was studied the problem on describing a class $\mathcal{X}$ of such Banach
spaces that $\left[  \mathsf{L}:(\mathcal{X}\subset X^{<f})\rightarrow
\mathcal{B}\right]  $, where $X^{f}$ denotes the class of all Banach spaces
that are finitely equivalent to a given Banach space $X$ (definitions see
below). This problem is completely similar to the classical one, which may be
formulated as a description of such class $\mathcal{X}$ of Banach spaces that
either $\left[  \mathsf{L}:(\mathcal{X}\subset\left(  l_{\infty}\right)
^{<f})\rightarrow\mathcal{B}\right]  $ or, equivalently, $\left[
\mathsf{L}:\left(  \mathcal{X}\subset\mathcal{B}\right)  \rightarrow
\mathcal{B}\right]  $ ($X^{<f}$ denotes the class of all Banach spaces that
are finitely representable in $X$).

The mentioned problem was partially solved in [6] where it was introduced a
class $\mathcal{E}(X^{f})$ that was called \textit{the class of existentially
closed spaces in the class }$X^{f}$; it was shown that for any
infinite-dimensional Banach space the class $\mathcal{E}(X^{f})$ is non-empty
(moreover, $\mathcal{E}(X^{f})$ is confinal with $X^{f}$) and, at last, that
for any superreflexive Banach space $X$ the following theorem is valid:
\[
\left[  \mathsf{L}:(\mathcal{E}(X^{f})\subset X^{f})\rightarrow\mathcal{B}%
\right]  .
\]

In the present paper it will be given the full description of the class
$\mathcal{V}$ that satisfy $\left[  \mathsf{L}:(\mathcal{V}\subset
X^{<f})\rightarrow\mathcal{B}\right]  $ for a given $X\in\mathcal{B}$.

Also, it will be considered the corresponding version of the classical problem
on extensions '\textit{to}', i.e., the problem of description of a such class
$\mathcal{W}\subset\mathcal{B}$ that $\left[  \mathsf{L}:(X^{<f}\subset
X^{<f})\rightarrow\mathcal{W}\right]  $.

The results, obtained in this paper and concerning the mentioned problem are
differ from classical on description a class $\mathcal{W}\subset\mathcal{B}$
such that $\left[  \mathsf{L}:\left(  \mathcal{B}\subset\mathcal{B}\right)
\rightarrow\mathcal{W}\right]  $ (which is, as it is known, equivalent to the
extensions 'from' problem). Namely, it will be shown that a desired class
$\mathcal{W}$ is non-empty only in special cases.

\section{Definitions and notations}

Let $\mathcal{B}$ denotes the class of all Banach spaces.

The \textit{dimension} $\dim(X)$ of $X\in\mathcal{B}$ is the least cardinality
of a subset $A\subset X$, which \textit{linear span} $\operatorname{lin}(A)$
is dense in $X$ (equivalently, of such $A\subset X$ that the \textit{closure}
$\overline{\operatorname{lin}(A)}$, which will be in the future denoted by
$\operatorname{span}(A)$, is the whole space $X$).

\begin{definition}
Let $X$, $Y$ are Banach spaces. $X$ is \textit{finitely representable} in $Y$
(in symbol: $X<_{f}Y$) if for each $\varepsilon>0$ and for every
finite-dimensional subspace $A$ of $X$ there exists a subspace $B$ of $Y$ and
an isomorphism $u:A\rightarrow B$ such that $\left\|  u\right\|  \left\|
u^{-1}\right\|  \leq1+\varepsilon$.

Spaces $X$ and$\ Y$ are said to be finitely equivalent, shortly: $X\sim_{f}Y$,
if $X<_{f}Y$ and $Y<_{f}X$.

Any Banach space $X$ generates a class
\[
X^{f}=\{Y\in\mathcal{B}:X\sim_{f}Y\}
\]
\end{definition}

Let $\mathcal{K}$ be a class of Banach spaces; $X$ be a Banach space. A notion
$X<_{f}\mathcal{K}$ (or, equivalently, $X\in\mathcal{K}^{<f}$ means that
$X<_{f}Y$ for every $Y\in\mathcal{K}$.

For any two Banach spaces $X$, $Y$ their \textit{Banach-Mazur distance }is
given by
\[
d(X,Y)=\inf\{\left\|  u\right\|  \left\|  u^{-1}\right\|  :u:X\rightarrow
Y\},
\]
where $u$ runs all isomorphisms between $X$ and $Y$ and is assumed, as usual,
that $\inf\varnothing=\infty$.

It is well known that $\log d(X,Y)$ defines a metric on each class of
isomorphic Banach spaces, where almost isometric Banach spaces are identified.

Recall that Banach spaces $X$ and $Y$ are said to be \textit{almost isometric}
provided $d(X,Y)=1$. Certainly, any pair of almost isometric
finite-dimensional Banach spaces $X$ and $Y$ is isometric.

The set $\frak{M}_{n}$ of all $n$-dimensional Banach spaces, equipped with
this metric, is a compact metric space, which is called \textit{the Minkowski
compact} $\frak{M}_{n}$.

A disjoint union $\cup\{\frak{M}_{n}:n<\infty\}=\frak{M}$ is a separable
metric space, which is called the \textit{Minkowski space}.

Consider a Banach space $X$. Let $H\left(  X\right)  $ be a set of all its
\textit{different} finite-dimensional subspaces (\textit{isometric
finite-dimensional subspaces of }$X$\textit{\ in }$H\left(  X\right)
$\textit{\ are identified}). Thus, $H\left(  X\right)  $ may be regarded as a
subset of $\frak{M}$, equipped with the restriction of the metric topology of
$\frak{M}$.

Of course, $H\left(  X\right)  $ need not to be a closed subset of $\frak{M}$.
Its closure in $\frak{M}$ will be denoted by $\overline{H\left(  X\right)  }$.

From the definitions it follows that $X<_{f}Y$ if and only if $\overline
{H\left(  X\right)  }\subseteq\overline{H\left(  Y\right)  }$.

Certainly, spaces $X$ and $Y$ are \textit{finitely equivalent }( $X\sim_{f}Y$)
if and only if $\overline{H\left(  X\right)  }=\overline{H\left(  Y\right)  }$.

Thus, there is a one to one correspondence between classes of finite
equivalence $X^{f}=\{Y\in\mathcal{B}:X\sim_{f}Y\}$ and closed subsets of
$\frak{M}$ of kind $\overline{H\left(  X\right)  }$.

Indeed, if $Y$ and $Z$ belong to $X^{f}$, then $S=\overline{H(Y)}%
=\overline{H(Z)}$. The set $S$, uniquely determined by $X$ (or, equivalently,
by $X^{f}$), will be denoted by $\frak{M}(X^{f})$ and will be referred to as
\textit{the Minkowski's base of the class} $X^{f}$. In what follows, for any
natural $n$ put $\frak{M}_{n}(X^{f})=\frak{M}(X^{f})\cap\frak{M}_{n}$.

\begin{definition}
Let $I$ be a set; $D$ be an ultrafilter over $I$; $\{X_{i}:i\in I\}$ be a
family of Banach spaces. An \textit{ultraproduct }$(X_{i})_{D}$ is given by a
quotient space
\[
(X)_{D}=l_{\infty}\left(  X_{i},I\right)  /N\left(  X_{i},D\right)  ,
\]
where $l_{\infty}\left(  X_{i},I\right)  $ is a Banach space of all families
$\frak{x}=\{x_{i}\in X_{i}:i\in I\}$ such that
\[
\left\|  \frak{x}\right\|  =\sup\{\left\|  x_{i}\right\|  _{X_{i}}:i\in
I\}<\infty;
\]
$N\left(  X_{i},D\right)  $ is a subspace of $l_{\infty}\left(  X_{i}%
,I\right)  $, which consists of such $\frak{x}$'s$\ $that
\[
\lim_{D}\left\|  x_{i}\right\|  _{X_{i}}=0.
\]
\end{definition}

If all $X_{i}$'s are equal to a space $X\in\mathcal{B}$ then the ultraproduct
is called the \textit{ultrapower} and is denoted by $\left(  X\right)  _{D}$.

The operator $d_{X}:X\rightarrow\left(  X\right)  _{D}$ that sends any $x\in
X$ to an element $\left(  x\right)  _{D}\in\left(  X\right)  _{D}$, which is
generated by a stationary family $\{x_{i}=x:i\in I\}$ is called the
\textit{canonical embedding }of $X$ into its ultrapower $\left(  X\right)
_{D}$.

A Banach space $X$ is finitely representable in a Banach space $Y$ if and only
if there exists such ultrafilter $D$ (over $I=\cup D$) that $X$ is isometric
to a subspace of the ultrapower $(Y)_{D}$ (cf. [7])

Let us note that for any class $W^{f}$ its Minkowski's base $\frak{N}%
=\frak{M}(W^{f})$ has following properties:

(\textbf{C}) $\frak{N}\ $\textit{is a closed subset of the Minkowski's space
}$\frak{M}$;

(\textbf{H}) \textit{If }$A\in\frak{N}$\textit{ and }$B\in H(A)$\textit{ then
}$B\in\frak{N}$;

(\textbf{A}$_{0}$)\textit{ For any }$A$\textit{, }$B\in\frak{N}$\textit{
\ there exists }$C\in\frak{N}$\textit{ \ such that }$A\in H(C)$\textit{ and}
$B\in H(C)$.

\begin{theorem}
Let $\frak{N}$ be a set of finite-dimensional Banach spaces; $\frak{N}%
\subset\frak{M}$. If $\frak{N}$ has properties (\textbf{C}), (\textbf{H}) and
(\textbf{A}$_{0}$) then there exists a class $X^{f}$ such that $\frak{N}%
=\frak{M}(X^{f})$.
\end{theorem}

\begin{proof}
Conditions (\textbf{H}) and (\textbf{A}$_{0}$) in a natural (not unique) way
define a partial order on $\frak{M}(X^{f})=\frak{N}$. Let $D$ be an
ultrafilter which is consistent with this order. Then the ultraproduct
$(\frak{N})_{D}$ of all spaces from $\frak{N}$ has desired properties. Since
the set $H\left(  (\frak{N})_{D}\right)  $ is closed, $\frak{M}(W^{f}%
)=\frak{N}$.
\end{proof}

The following classifications of Banach spaces due to L. Schwartz [8].

For a Banach space $X$ its $l_{p}$-\textit{spectrum }$S(X)$ is given by
\[
S(X)=\{p\in\lbrack0,\infty]:l_{p}<_{f}X\}.
\]

Certainly, if $X\sim_{f}Y$ then $S(X)=S(Y)$. Thus, the $l_{p}$-spectrum $S(X)
$ may be regarded as a property of the whole class $X^{f}$. So, notations like
$S(X^{f})$ are of obvious meaning.

Let $X$ be a Banach space. It is said to be:

\begin{itemize}
\item $c$-\textit{convex,} if $\infty\notin S(X)$;

\item $B$-\textit{convex,} if $1\notin S\left(  X\right)  $;

\item \textit{finitely universal,} if $\infty\in S(X)$.
\end{itemize}

The $l_{p}$-spectrum is closely connected with notions of type and cotype.
Recall the definition.

\begin{definition}
Let $1\leq p\leq2\leq q\leq\infty$. A Banach space $X$ is said to be of
\textit{type} $p$, respectively, of \textit{cotype }$q$, if for every finite
sequence $\{x_{n}:n<N\}$ of its elements
\[
\int_{0}^{1}\left\|  \sum\nolimits_{n=0}^{N-1}r_{n}\left(  t\right)
x_{n}\right\|  dt\leq t_{p}\left(  X\right)  \left(  \sum\nolimits_{n=0}%
^{N-1}\left\|  x_{n}\right\|  ^{p}\right)  ^{1/p},
\]
respectively,%
\[
\left(  \sum\nolimits_{n=0}^{N-1}\left\|  x_{n}\right\|  ^{q}\right)
^{1/q}\leq c_{q}\left(  X\right)  \int_{0}^{1}\left\|  \sum\nolimits_{n=0}%
^{N-1}r_{n}\left(  t\right)  x_{n}\right\|  dt,
\]
where $\left\{  r_{n}\left(  t\right)  :n<\infty\right\}  $ are Rademacher functions.
\end{definition}

When $q=\infty$, the sum $\left(  \sum_{n=0}^{N-1}\left\|  x_{n}\right\|
^{q}\right)  ^{1/q}$ is replaced with $\sup\nolimits_{n<N}\left\|
x_{n}\right\|  $. Constants $t_{p}\left(  X\right)  $ and $c_{q}\left(
X\right)  $\ in these inequalities depend only on $X$. Their least values
$T_{p}\left(  X\right)  =\inf t_{p}\left(  X\right)  $ and $C_{q}\left(
X\right)  =\inf c_{q}\left(  X\right)  $, are called \textit{the type}
$p${\small -}\textit{constant} $T_{p}\left(  X\right)  $ and \textit{the
cotype} $q$\textit{-constant }$C_{q}\left(  X\right)  $.

Every Banach space is of type $1$ and of cotype $\infty$.

If $X$ is of type $p$ and of cotype $q$ with the constants $T_{p}\left(
X\right)  =T$, $C_{q}\left(  X\right)  =C$ than any $Y\in X^{f}$ is of same
type and cotype and its type-cotype constants are equal to those of $X$. Thus,
it may be spoken about the type and cotype of the whole class $X^{f}$ of
finite equivalence. Notice that $X$ and $Y$ are isomorphic then these spaces
are of the same type and cotype.

It is known (see [5]) that
\begin{align*}
\inf S(X)  &  =\sup\{p\in\lbrack1,2]:T_{p}\left(  X\right)  <\infty\};\\
\sup S(X)  &  =\inf\{q\in\lbrack2,\infty]:C_{q}\left(  X\right)  <\infty\}.
\end{align*}

\begin{definition}
A Banach space $X$ is said to be \textit{superreflexive} if every space of the
class $X^{f}$ is reflexive.
\end{definition}

Equivalently, $X$ is superreflexive if any $Y<_{f}X$ is reflexive. Clearly,
any superreflexive Banach space is $B$-convex. Since $S(X)$ is either
$[1,\infty]$ or is a closed subset of $[1,\infty)$ (cf. [6]), any
superreflexive space is of \textit{non-trivial} (i.e., non equal to $1$) type.

Below it will be used the notion of direct (inductive) limit of Banach spaces.

\begin{definition}
Let $\left\langle I,\ll\right\rangle $ be a \textit{partially ordered set}. It
said to be \textit{directed} (to the right hand) if for any $i,j\in I$ there
exists $k\in I$ such that $i\ll k$ and $j\ll k$ .
\end{definition}

Let $\left\{  X_{i}:i\in I\right\}  $ be a set of Banach spaces that are
indexed by elements of an directed set $\left\langle I,\ll\right\rangle $ (any
set of elements of arbitrary nature, which is indexed by elements of the
directed set will be called a \textit{direction}). Let $m_{i,j}:X_{i}%
\rightarrow X_{j}$ $(i\ll j)$ be isomorphic embeddings.

\begin{definition}
A system $\left\{  X_{i},m_{i,j}:i,j\in I;i\ll j\right\}  $ is said to be an
\textit{inductive (}or\textit{\ direct) system} if
\[
m_{i,i}=Id_{X_{i}};\text{ \ \ }m_{i,k}=m_{j,k}\cdot m_{i,j}%
\]
for all $i\ll j\ll k$ ($Id_{Y}$ denotes the identical operator on $Y$).
\end{definition}

Its \textit{inductive }(or \textit{direct}) \textit{limit}, $\underset
{\rightarrow}{\lim}X_{i}$ is defined as follows. Put
\[
X=\cup\left\{  X_{i}\times\left\{  i\right\}  :i\in I\right\}
\]

Elements of $X$ are pairs $(x,i),$ where $x\in X_{i}$. Let $=_{eq}$ be a
relation of equivalence of elements of $X$, which is given by the following rule:

\begin{center}
$(x,i)=_{eq}(y,j)$ if $m_{i,k}x=m_{j,k}y$ for some $k\in I$.
\end{center}

A class of all elements of $X$ that are equivalent to a given $(x,i)$ is
denoted by
\[
\lbrack x,i]=\{(y,j):(y,j)=_{eq}(x,i)\}.
\]

A set of all equivalence classes $[x,i]$ is denoted $X_{\infty}$. Clearly,
$X_{\infty}$ is a linear space. Let $\left\|  [x,i]\right\|  =\lim_{I}\left\|
m_{i,j}x\right\|  _{X_{j}}$ be a semi-norm on $X_{\infty}$ ($\lim_{I}$ denotes
the limit by the directed set $I$). Put
\[
Null(X)=\{[x,i]:\left\|  [x,i]\right\|  =0\}.
\]

\begin{definition}
The direct limit of the inductive system $\left\{  X_{i},m_{i,j}:i\ll j\in
I\right\}  $ is a quotient space
\[
\underset{\rightarrow}{\lim}X_{i}=X_{\infty}/Null(X).
\]
\end{definition}

Clearly, if all spaces $X_{i}$ ($i\in I$) belong to a given class $Z^{f}$,
then $\underset{\rightarrow}{\lim}X_{i}\in Z^{f}$; if all $X_{i}$'s are
finitely representable in $Z$ then $\underset{\rightarrow}{\lim}X_{i}<_{f}Z$ too.

A particular case of direct limits is the case when all embeddings. $\left(
m_{i,j}\right)  $ are isometric and all spaces $\left(  X_{i}\right)  $ from
the inductive system are directed into the chain:
\[
X_{0}\hookrightarrow X_{1}\hookrightarrow...\hookrightarrow X_{i}%
\hookrightarrow...\text{ .}%
\]

In this case $\underset{\rightarrow}{\lim}X_{i}$ may be identified with the
closure of the union of this chain:%
\[
\underset{\rightarrow}{\lim}X_{i}=\overline{\cup_{i<\infty}X_{i}}.
\]

\section{Algebraically closed Banach spaces}

Recall that for a Banach space $X$ and its subspace $A$, the
\textit{annihilator} $A^{\perp}$ of $A$ in $X^{\ast}$ is given by:%
\[
A^{\perp}=\{f\in B^{\ast}:f\left(  a\right)  =0\text{ \ for all \ }a\in A\}.
\]

\begin{definition}
$A$ is said to be an \textit{algebraic subspace} of $X$ if there is a
projection of norm one from $X^{\ast}$ onto the annihilator $A^{\perp}$.
\end{definition}

An equivalent definition of algebraic subspaces according to [9] looks like:

\begin{itemize}
\item $A$\textit{ is an algebraic subspace of }$X$\textit{ if for every
}$\varepsilon>0$\textit{ and each finite-dimensional subspace }$B$\textit{ of
}$X$\textit{ there exists a surjection }$u:B\rightarrow A$\textit{ of norm
}$\left\|  u\right\|  \leq1+\varepsilon$\textit{ which is identical on the
intersection }$B\cap A$\textit{: }$u\mid_{B\cap A}=Id_{B\cap A}$\textit{.}
\end{itemize}

\begin{definition}
A Banach space $E$ is said to be algebraically closed in the class $X^{f}$ if
for any isometric embedding $i:E\rightarrow Z$ into an arbitrary space $Z\in
X^{f}$ its image $iE$ is an algebraic subspace of $Z$, or, equivalently, a
double adjoint operator $i^{\ast\ast}:E^{\ast\ast}\rightarrow Z^{\ast\ast}$
maps $E^{\ast\ast}$ onto a norm one complemented subspace of $Z^{\ast\ast}$:
there exists a projection $P:Z^{\ast\ast}\rightarrow E^{\ast\ast}$ of norm
$\left\|  P\right\|  =1$.
\end{definition}

The class of all spaces $E$, algebraically closed in $X^{f}$, is denoted by
$\mathcal{A}(X^{f})$.

\begin{definition}
($[9]$). Let $X$ be a Banach space; $Y$ be its subspace. $Y$ is said to be a
reflecting subspace of $X$ (symbolically: $Y\prec_{u}X$) if for every
$\varepsilon>0$ and every finite-dimensional subspace $A\hookrightarrow X$
there exists an isomorphic embedding $u:A\rightarrow Y$ such that $\left\|
u\right\|  \left\|  u^{-1}\right\|  \leq1+\varepsilon$ and $u\mid_{A\cap
Y}=Id_{A\cap Y}$.
\end{definition}

As it was shown in [9], if $Y\prec_{u}X$ then $Y^{\ast\ast}$ is an image of a
norm one projection $P:X^{\ast\ast}\rightarrow Y^{\ast\ast}$ (under canonical
embedding of $Y^{\ast\ast}\ $into $X^{\ast\ast}$).

\begin{definition}
($[6]$). A Banach space $E$ is said to be existentially closed in the class
$X^{f}$ if for any isometric embedding $i:E\rightarrow Z$ into an arbitrary
space $Z\in X^{f}$ its image $iE$ is a reflecting subspace of $Z$:
$iY\prec_{u}Z$.
\end{definition}

The class of all spaces $E$, existentially closed in $X^{f}$, is denoted by
$\mathcal{E}(X^{f})$.

In [6] it was shown that for any Banach space $X$ the class $\mathcal{E}%
(X^{f})$ is non-empty; moreover, any $Y<_{f}X^{f}$ may be isometrically
embedded into some $E\in\mathcal{E}(X^{f})$ of dimension $\dim(E)=\max
\{\dim(Y),\omega\}$ ($\omega$ denotes the first infinite ordinal number). In
fact, the existence of sufficiently many existentially closed spaces in every
class $\mathcal{K}$ of Banach spaces, which is closed under unions of chains
(and their closures) is an easy consequence of the well known Lindstr\"{o}m
theorem from the model theory.

However, to apply this result to the Banach space theory it will be needed
some additional constructions because of axioms of Banach spaces (besides the
triangle inequality) cannot be represented in the first order logic.

The proof, which is presented below completely eliminates any metamathematics
although it based on some model-theoretical ideas of M. Yasuhara [10].

\begin{theorem}
For every infinite-dimensional Banach space $X$ the class $\mathcal{E}(X^{f})$
of all existentially closed in $X^{f}$ Banach spaces is non-empty. Every $Y\in
X^{f}$ may be isometrically embedded into some $E\in\mathcal{E}(X^{f})$ of the
dimension $\dim(E)=\dim(Y)$.
\end{theorem}

\begin{proof}
A triple $\left(  u,A,Z\right)  $ will be called an $\varepsilon$-triple over
$Y$ \ if $A\hookrightarrow Y$ is a finite-dimensional subspace; $Z\in
\frak{M}(X^{f})$ and $u:A\rightarrow Z$ be an $\left(  1+\varepsilon\right)
$-isomorphic embedding (i.e., $\left\|  u\right\|  \left\|  u^{-1}\right\|
\leq1+\varepsilon$).

It will be said that the $\varepsilon$-triple $\left(  u,A,Z\right)  $ is
realized by $W\in X^{f}$, which contains $Y$, if there exists such isometric
embedding $v:Z\rightarrow W$ that for all elements $y\in A$ the composition
$v\circ u$ fixes $y$: $\left(  v\circ u\right)  y=y$. Surely, if $\left(
u,A,Z\right)  $ is realized by $W$, it is realized by every space $W_{0}$ that
contains $\operatorname{span}\{Y\cup vZ\}$, in particular, there exists a such
$W_{0}$, which dimension is equal to $\dim Y.$ In what follows it will be
constantly assuming that the dimension of all realizing spaces $W$ in the
proof satisfies $\dim W=\dim Y=\varkappa$ for some cardinal $\varkappa$.

Chose a dense subset $\left(  y_{\gamma}\right)  _{\gamma<\varkappa}$ of the
unit sphere $\partial B\left(  Y\right)  =\{y\in Y:\left\|  y\right\|  =1\}$;
a sequence $\left(  \varepsilon_{n}\right)  _{n<\infty}$ with $\varepsilon
_{n}\searrow0$ and number by the cardinal $\varkappa$ all $\varepsilon_{n}%
$-triples $\left(  u,A_{I},Z\right)  $ over $Y$, $n<\infty$, where $A_{I}$ is
spanned by a finite subset of $\left(  y_{\gamma}\right)  _{\gamma<\varkappa}%
$: $A=\operatorname{span}\{\left(  y_{\gamma}\right)  _{\gamma<\varkappa
}:I\subset\varkappa\}$ and $Z$ is chosen from the subset $\frak{N}\subset
$.$\frak{M}(X^{f})$ that consists only of spaces, which are isometric to
spaces $Y_{I}$ of kind
\[
Y_{I}=\operatorname{span}\{\left(  y_{\gamma}\right)  _{\gamma<\varkappa
}:I\subset\varkappa;\operatorname{card}I<\infty\}.
\]

Certainly, $\frak{N}\ $\ is a dense subset of $\frak{M}(X^{f})$.

So, let $\frak{A}$ $=\{\left(  u_{\gamma},A_{\gamma},Z_{\gamma}\right)
:\gamma<\varkappa\}$ be a numeration of all triples over $Y$ of the described kind.

Define two sequences of Banach spaces: $\left(  Y_{\gamma}\right)
_{\gamma<\varkappa}$ and $\left(  \overline{Y}_{\gamma}\right)  _{\gamma
<\varkappa}$ inductively.

Put $\overline{Y}_{0}=Y$.

If $\left(  Y_{\gamma}\right)  _{\gamma<\delta(<\varkappa)}$ are already
defined, put
\[
\overline{Y}_{\delta}=\overline{\cup\{Y_{\gamma}:\gamma<\delta\}.}%
\]

Let $\left(  u_{\delta},A_{\delta},Z_{\delta}\right)  \in\frak{A}$ be the
$\delta$'th triple over $Y$. It may be assumed to be an $\varepsilon_{n\left(
\delta\right)  }$-triple over $\overline{Y}_{\delta}$. If this triple is
realized by some $W$, put $Y_{\delta}=W$. In a contrary case put $Y_{\delta
}=\overline{Y}_{\delta}$.

Put%
\begin{align*}
Y^{\left(  1\right)  }  &  =\underset{\rightarrow}{\lim}Y_{\gamma};\\
Y^{\left(  n+1\right)  }  &  =(Y^{\left(  n\right)  })^{\left(  1\right)  };\\
Y_{\infty}  &  =\underset{\rightarrow}{\lim}Y^{\left(  n\right)  }.
\end{align*}

Certainly, $Y\hookrightarrow Y_{\infty}$; $\dim Y_{\infty}=\varkappa$ and
$Y_{\infty}\in Y^{f}=X^{f}$. Let us show that $Y_{\infty}\in\mathcal{E}%
(X^{f})$.

Assume that $V\in X^{f}$ and $Y_{\infty}\hookrightarrow V$.

Chose a finite-dimensional subspace $Z\hookrightarrow V$. Let
$A\hookrightarrow Y_{\infty}$ be chosen from the collection of spaces of
$\frak{A}$. Let $u:A\rightarrow Y_{\infty}$ be the identical embedding and
$B=\operatorname{span}\left(  A\cup Z\right)  $. Since $\frak{N}\ $\ is a
dense subset of $\frak{M}\left(  X^{f}\right)  $, for every $\varepsilon>0$
the triple $\left(  u,A,B\right)  $ may be regarded as an $\varepsilon$-triple
over $Y_{\infty}$, which is realized by $V$. Certainly, it may be assumed that
this triple is a triple over $Y^{\left(  n\right)  }$ for some $n<\infty$ and,
hence, we meet it at some step $\delta<\varkappa$ of our construction as
$\left(  u_{\delta},A_{\delta},Z_{\delta}\right)  $. On the $\delta$'th step
of construction of $Y^{\left(  n+1\right)  }=\left(  Y^{\left(  n\right)
}\right)  ^{\left(  1\right)  }$ it was assumed to be a triple over
$\overline{(Y^{\left(  n\right)  })}_{\delta}$. Since this triple is realized
by $V$, then the space $\left(  Y^{\left(  n\right)  }\right)  _{\delta}$ also
realizes it.

Hence, for every $\varepsilon>0$ there exists a such embedding $v:B\rightarrow
(Y^{\left(  n\left(  \varepsilon\right)  \right)  })_{\delta\left(
\varepsilon\right)  }$ that $\left(  v\circ u\right)  x=x$ for all $x\in A$;
$\left\|  v\right\|  \left\|  v^{-1}\right\|  \leq1+\varepsilon$.

Since $Z\hookrightarrow V$ is arbitrary, this implies that $Y_{\infty}$ is a
reflecting subspace of $V$. Because of $V\in X^{f}$ is arbitrary too,
$Y_{\infty}\in\mathcal{E}(X^{f})$.
\end{proof}

It follows from definitions that $\mathcal{E}(X^{f})\subseteq\mathcal{A}%
(X^{f})$. So, $\mathcal{A}(X^{f})$ also is non-empty and confinal with $X^{f}$
for every Banach space $X$. Clearly, $\mathcal{E}(X^{f})\subseteq X^{f}$.

Put $\mathcal{A}_{0}(X^{f})=\mathcal{A}(X^{f})\cap X^{f}$.

The following inclusions are obvious:
\[
\mathcal{E}\left(  X^{f}\right)  \subseteq\mathcal{A}_{0}\left(  X^{f}\right)
\subseteq\mathcal{A}\left(  X^{f}\right)  .
\]

\section{Universally complemented Banach spaces}

\begin{definition}
A Banach space $E$ is said to be universally complemented in the class $X^{f}$
if for any isometric embedding $i:E\rightarrow Z$ into an arbitrary space
$Z\in X^{f}$ its image $iE$ admits a norm one projection $P:Z\rightarrow iY$.
\end{definition}

The class of all spaces $E$ that are universally complemented in a class
$X^{f}$ is denoted by $\frak{P}_{1}(X^{f})$.

Certainly, $E\in\frak{P}_{1}(X^{f})$ if and only if for any isometric
embedding $i:E\rightarrow Z$ into an arbitrary space $Z\in X^{f}$ and for any
Banach space $W$ every operator $u:E\rightarrow W$ may be extended to an
operator $V:Z\rightarrow W$ of the same norm, i.e., there exists such
$V:Z\rightarrow W$ that $V\mid_{iY}=u$ and $\left\|  V\right\|  =\left\|
u\right\|  $.

In other words, $\frak{P}_{1}(X^{f})$ is the maximal class, which satisfies
the condition
\[
\left[  \mathsf{L}:(\frak{P}_{1}(X^{f})\subset X^{f})\rightarrow
\mathcal{B}\right]  .
\]

Below it will be also considered the class $\frak{P}_{1}^{\left(  c\right)
}\left(  X^{f}\right)  $, which is given by
\[
\left[  \mathsf{L}^{\left(  c\right)  }:(\frak{P}_{1}^{\left(  c\right)
}(X^{f})\subset X^{f})\rightarrow\mathcal{B}\right]
\]
and $\frak{P}_{1}^{wc}(X^{f})$ such that
\[
\left[  \mathsf{L}^{\left(  wc\right)  }:(\frak{P}_{1}^{\left(  wc\right)
}(X^{f})\subset X^{f})\rightarrow\mathcal{B}\right]  ,
\]
where $\mathsf{L}^{\left(  c\right)  }$ (resp., $\mathsf{L}^{\left(
wc\right)  }$) denotes the class of all compact (resp., weakly compact) operators.

As it was shown in [6], for any Banach space $X$ the following set-theoretical
inclusions are true:
\[
\text{If }\mathit{X}\text{ is superreflexive then }\mathcal{E}(X^{f}%
)\subseteq\frak{P}_{1}(X^{f});
\]%
\[
\text{In a general case }\mathcal{E}(X^{f})\subseteq\frak{P}_{1}^{\left(
c\right)  }(X^{f})\subseteq\frak{P}_{1}^{\left(  wc\right)  }(X^{f}).
\]
Below it will be proved a more exact result. Its formulation needs to
introduce one more notion.

\begin{definition}
A Banach space $X$ is said to be quasi-dual if there exists a projection of
norm one from $X^{\ast\ast}$ onto $k_{X}X$, where $k_{X}$ denotes the
canonical embedding of $X$ into its second conjugate.
\end{definition}

A class of all quasi-dual Banach spaces will be denoted by $\frak{qd}$. The
following proposition belongs to the folklore.

\begin{proposition}
$X\in\frak{qd}$ if and only if there exists a conjugate Banach space $Y^{\ast
}$, an isometric embedding $i:X\rightarrow Y^{\ast}$ and a projection
$P:Y^{\ast}\rightarrow iX$ of norm $\left\|  P\right\|  =1$.
\end{proposition}

\begin{proof}
If $X\in\frak{qd}$, let $Y=X^{\ast}$. Conversely, let $P:Y^{\ast}\rightarrow
iX$ \ be a projection of norm one. Put $R=P\circ Q:Y^{\ast\ast\ast}\rightarrow
X$, where $Q:Y^{\ast\ast\ast}\rightarrow k_{Y^{\ast}}Y^{\ast}$ is a canonical
projection. Certainly, $\left\|  R\right\|  =1$ and the restriction of $R$ to
$X^{\perp\perp}=X^{\ast\ast}$ is the desired projection from $X^{\ast\ast}$
onto $X$.
\end{proof}

The main reason in introducing algebraically closed Banach spaces describes
the following result.

\begin{theorem}
$\frak{P}_{1}(X^{f})=\mathcal{A}(X^{f})\cap\frak{qd}$.
\end{theorem}

\begin{proof}
Let $Y\in\mathcal{A}(X^{f})\cap\frak{qd}$; $W\in X^{f}$; $j:Y\rightarrow W$ be
an isometric embedding. Then there exists a projection $P:W^{\ast\ast
}\rightarrow j^{\ast\ast}Y^{\ast\ast}$ of norm one and a projection
$Q:Y^{\ast\ast}\rightarrow k_{Y}Y$. Certainly, $R=\left(  k_{W}\right)
^{-1}\circ P\circ Q$ is a projection from $W$ onto $jY$.

Conversely, let $Y\in\frak{P}_{1}(X^{f})$, $W\in X^{f}$; $j:Y\rightarrow W$ be
an isometric embedding. Then there is a projection $R:W\rightarrow jY$. Since
$W$ is arbitrary, $Y\in\mathcal{A}(X^{f})$. Because of $Y^{\ast\ast}\in X^{f}%
$, $Y\in\frak{qd}$. Hence, $Y\in\mathcal{A}(X^{f})\cap\frak{qd}$.
\end{proof}

Following [5] introduce one more notion.

\begin{definition}
Let $X$ be a Banach space; $Z<_{f}X$. It will be said that $Z$ has the
identity extension property in the class $X^{f}$, shortly: $Z\in IEP(X^{f})$,
if for any $W\in X^{f}$ that contains $Z$ as a subspace, the identical
operator $Id_{Z}:Z\rightarrow Z$ may be extended to the norm one operator
$T:W\rightarrow Z^{\ast\ast}$.
\end{definition}

\begin{theorem}
If $Z\in\mathcal{A}(X^{f})$ then $Z\in IEP(X^{f})$.
\end{theorem}

\begin{proof}
Let $W\in X^{f}$; $Z\in\mathcal{A}(X^{f})$. Consider a set $\frak{G}$ of pairs
$\left(  B,\varepsilon\right)  $, where $B$ runs a set $G(W)$ of all
finite-dimensional subspaces of $W$ (isometric subspaces in $G(W)$ are not
identified !); $\varepsilon>0$. Define on $\frak{G}$ a partial ordering $\ll$,
assuming that $\left(  B_{1},\varepsilon_{1}\right)  \ll\left(  B_{2}%
,\varepsilon_{2}\right)  $ if $B_{1}\hookrightarrow B_{2}$ and $\varepsilon
_{2}\leq\varepsilon_{1}$.

Since $Z\in\mathcal{A}(X^{f})$, there exists a projection $P:W^{\ast\ast
}\rightarrow Z^{\ast\ast}$ ($Z^{\ast\ast}$ and $Z^{\perp\perp}$ are
identified) of norm one; its restriction to $B\longrightarrow W$ maps $B$ into
$Z^{\ast\ast}$.

Certainly, $P\mid_{B\cap Z^{\ast\ast}}=Id_{B\cap Z^{\ast\ast}}$.

Since $B\cap Z^{\ast\ast}\hookrightarrow Z^{\ast\ast}$, by the principle of
local reflexivity [11] for every $\varepsilon>0$ there exists an operator
$T=T_{B,\varepsilon}:B\cap Z^{\ast\ast}\rightarrow Z$ such that $\left\|
T\right\|  \left\|  T^{-1}\right\|  \leq1+\varepsilon$ and $T\mid_{B\cap
Z}=Id\mid_{B\cap Z}$. Notice that there exists a projection $Q=P_{B}%
:B\rightarrow B\cap Z^{\ast}$ of norm one. Put $\widetilde{T}_{B,\varepsilon
}=T_{B,\varepsilon}\circ Q$. Immediately $\widetilde{T}_{B,\varepsilon
}:B\rightarrow Z$; $\left\|  \widetilde{T}_{B,\varepsilon}\right\|  =\left\|
T_{B,\varepsilon}\right\|  $.

Let $r>0$ be a real number. Put
\[
U(r)=\{z\in Z^{\ast\ast}:\left\|  z\right\|  =r\}.
\]
This set is weakly* compact.

Hence, by the Tychonoff theorem, a product $\Pi=\prod\nolimits_{w\in
W}U\left(  2\left\|  w\right\|  \right)  $ is compact in a product topology.

Let $t(w)$ denotes the $w$'s coordinate of a point $t\in\Pi$. Every operator
$\widetilde{T}_{B,\varepsilon}$ (which may be regarded as an operator from $B$
to $Z^{\ast\ast}$) defines a point $t_{B,\varepsilon}\in\Pi$, which is given
by%
\begin{align*}
t_{B,\varepsilon}\left(  w\right)   &  =\widetilde{T}_{B,\varepsilon}\left(
w\right)  \text{ \ for }w\in B;\\
&  =0\text{ \ \ \ \ \ \ \ \ \ \ \ for }w\notin B.
\end{align*}

Notice that $\left\langle \frak{G},\ll\right\rangle $ is the directed set.

Let $t$ be a limit point of a direction $t_{B,\varepsilon}$. It has properties:

\begin{itemize}
\item $t\left(  z\right)  =z$ \ for $z\in Z$.

\item $t\left(  aw_{1}+bw_{2}\right)  =at\left(  w_{1}\right)  +bt\left(
w_{2}\right)  $ for $a$, $b\in\mathbb{R}$; $w_{1}$, $w_{2}\in W$;

\item $\left\|  t\left(  w\right)  \right\|  \leq\left\|  w\right\|  $ for all
$w\in W$.
\end{itemize}

Let $T_{0}:W\rightarrow Z^{\ast\ast}$ be given by $T_{0}\left(  w\right)
=t\left(  w\right)  $ for all $w\in W$. Clearly, $T_{0}$ is the desired
extension of $Id_{Z}$.
\end{proof}

\begin{theorem}
Let $X$ be a Banach space, $Z<_{f}X$ and $Z\in IEP(X^{f})$. Then $Z^{\ast\ast
}\in\frak{P}_{1}(X^{f})$.
\end{theorem}

\begin{proof}
Let $Y\in X^{f}$; $j:Z^{\ast\ast}\rightarrow Y$ be an isometric embedding. For
simplicity $Z^{\ast\ast}$ will be regarded as a subspace of $Y$. Let
$W\in\mathcal{A}(X^{f})$ contains an isometric image of $Y$. As $W$ it may be
chosen the corresponding space from $\mathcal{E}(X^{f})$. For simplicity, $Y$
is assumed to be a subspace of $W$.

So, we have a chain $Z\hookrightarrow Z^{\ast\ast}\hookrightarrow
Y\hookrightarrow W$.

Since $Z\in IEP(X^{f})$, there exists an operator $T:W\rightarrow W$,
$TW\hookrightarrow Z^{\ast\ast}$ such that $T\mid_{Z}=Id_{Z}$; $\left\|
T\right\|  =1$. Clearly,%
\[
T^{\ast\ast}\left(  W^{\ast\ast}\right)  \hookrightarrow Z^{\ast\ast\perp
\perp};\text{ \ }T^{\ast\ast}\mid_{Z^{\ast\ast}}=Id_{Z^{\ast\ast}}.
\]

Let $Q:Z^{\ast\ast\ast}\rightarrow Z^{\ast}$ be a natural projection.

Then $Q^{-1}\left(  0\right)  =Z^{\perp}$ and $Q^{\ast}:Z^{\ast\ast\ast\ast
}\rightarrow Z^{\perp\perp}$ is also a projection.\ Since $Z^{\ast\ast
\perp\perp}$ and $Z^{\ast\ast\ast\ast}$ are isometric, there exists a
projection $P:Z^{\ast\ast\perp\perp}\rightarrow Z^{\ast\ast}$ of norm
$\left\|  P\right\|  =1$.

Certainly, the composition $P\circ T^{\ast\ast}:W^{\ast\ast}\rightarrow
Z^{\ast\ast}$ is a projection of norm one. Its restriction $R=\left(  P\circ
T^{\ast\ast}\right)  \mid_{Y}$ is a projection from $Y$ onto $Z^{\ast\ast}$;
$\left\|  R\right\|  =1$.

Since $Y\in X^{f}$ and $j:Z^{\ast\ast}\rightarrow Y$ are arbitrary,
$Z^{\ast\ast}\in\frak{P}_{1}(X^{f})$.
\end{proof}

\begin{corollary}
For any $X\in\mathcal{B}$ the class $\frak{P}_{1}(X^{f})$ is non-empty. If $X$
is superreflexive, it is identical with $\mathcal{A}(X^{f})$ (that contains
$\mathcal{E}(X^{f})$). If $X$ is not superreflexive and $c$-convex then it is
equal to $\mathcal{A}(X^{f})\cap\frak{qd}$ and contains all spaces of kind
$E^{\ast\ast}$, where $E\in\mathcal{A}(X^{f})$. If $X$ is not $c$-convex, then
$\frak{P}_{1}(X^{f})$ is exactly the class of all quasi-dual $\mathcal{L}%
_{\infty,1+0}$-spaces.
\end{corollary}

\begin{proof}
Is a consequence of the previous theorem. If $X$ is not $c$-convex, this is
just the classical result from [11].
\end{proof}

\begin{remark}
By the way, it was proved that for every space $E\in\mathcal{A}(X^{f})$ its
second conjugate $E^{\ast\ast}$ also belongs to $\mathcal{A}(X^{f})$.

Notice that the\ validity of the implication $E\in\mathcal{E}(X^{f}%
)\Rightarrow E^{\ast\ast}\in\mathcal{E}(X^{f})$ is doubtful.
\end{remark}

Certainly, in a non-reflexive case, not every space $E\in\mathcal{A}(X^{f})$
belongs to $\frak{P}_{1}(X^{f})$. However, every $E\in\mathcal{A}(X^{f})$ has
a weaker property of extension of compact and weakly compact operators.

\begin{theorem}
For any Banach space $X$ the class of all algebraically closed spaces
$\mathcal{A}(X^{f})$ is identical with $\frak{P}_{1}^{\left(  wc\right)
}(X^{f})$ (and also with $\frak{P}_{1}^{\left(  c\right)  }(X^{f})$).
\end{theorem}

\begin{proof}
Let $Y\in\mathcal{A}(X^{f})$; $Y\hookrightarrow Z\in X^{f}$; $u:Y\rightarrow
W$ be a weakly compact operator from $Y$ to an arbitrary Banach space $W$.
Hence, $u$ may be factored through a reflexive Banach space, say $V$ [12]:
$u=v\circ w$, where $w:Y\rightarrow V$, $v:V\rightarrow W$. Certainly,
$u^{\ast\ast}$ admits a factorization through the same space $V$; $u^{\ast
\ast}=v^{\ast\ast}\circ w^{\ast\ast}$; $w^{\ast\ast}:Y^{\ast\ast}\rightarrow
V$; $v^{\ast\ast}:V\rightarrow W^{\ast\ast}$. Since $Y\in\mathcal{A}(X^{f})$,
$Y$ has the property $IEP(X^{f})$. Let $u_{0}:Y\rightarrow Y^{\ast\ast
}\hookrightarrow Z^{\ast\ast}$ be the corresponding extension of identity
$Id_{Z}$.

Clearly, $v\circ w^{\ast\ast}\circ u_{0}$ is the desired extension of $u$.

Let $Y\in\frak{P}_{1}^{\left(  wc\right)  }(X^{f})$. Then $Y^{\ast\ast}%
\in\frak{P}_{1}(X^{f})$ and, hence, $Y^{\ast\ast}\in\mathcal{A}(X^{f})$. The
inclusion $Y\in\mathcal{A}(X^{f})$ follows from the principle of local reflexivity.
\end{proof}

\begin{summary}
For every Banach space $X$%
\[
\mathcal{A}(X^{f})\cap\frak{qd}=\frak{P}_{1}(X^{f})\text{; \ }\mathcal{A}%
(X^{f})=\frak{P}_{1}^{\left(  wc\right)  }(X^{f})=\frak{P}_{1}^{\left(
c\right)  }(X^{f}).
\]
The class $\mathcal{A}(X^{f})\cap\frak{qd}$ is not empty because of
$\mathcal{E}(X^{f})\subseteq\mathcal{A}(X^{f})$ and $E\in\mathcal{E}(X^{f}%
)\ $implies that $E^{\ast\ast}\in\mathcal{A}(X^{f})\cap\frak{qd}$.
\end{summary}

\section{Divisible classes of finite equivalence}

This section as well as the next one does not touch problems on extension.
Nevertheless their results are of great significance for these problems, as it
will be shown below.

Let $X$, $Y$ be Banach spaces. Recall that their \textit{algebraic direct sum}
is a vector space $X\dotplus Y$ of all pairs $\left(  x,y\right)  $ where
$x\in X$, $y\in Y$.

A norm $N(x,y)$ on this space, which turns it into a Banach space must satisfy
the following condition:

If $\underset{n\rightarrow\infty}{\lim}\left\|  x-x_{n}\right\|  _{X}=0$ and
$\underset{m\rightarrow\infty}{\lim}\left\|  y-y_{m}\right\|  _{Y}=0$ then%
\[
\underset{n,m\rightarrow\infty}{\lim}N(x-x_{n},y-y_{m})=0.
\]

The space $X\dotplus Y$ equipped with such norm $N$ will be called the
\textit{direct }$N$\textit{-sum} (or, simply the \textit{direct sum}) of $X$
and $Y$ and will be denoted by $X\oplus_{N}Y$. Without loss of generality it
will be assumed that every norm $N\left(  \cdot,\cdot\right)  $ under
consideration satisfies the following conditions:%
\[
N\left(  x,0\right)  =\left\|  x\right\|  _{X};\text{ \ }N\left(  0,y\right)
=\left\|  y\right\|  _{Y}.
\]

Operators\ $i:X\rightarrow X\oplus_{N}Y$ and $j:Y\rightarrow X\oplus_{N}Y$,
which are given by
\[
ix=\left(  x,0\right)  ,\text{ \ }x\in X;\text{ \ }jy=\left(  0,y\right)
,\text{ y}\in Y
\]
are called \textit{canonical }(isometric)\textit{ embeddings.}

\textit{Canonical projections }$P:$ $X\oplus_{N}Y\rightarrow iX$ and
$Q:X\oplus_{N}Y\rightarrow jY$ are given by%
\[
P\left(  x,y\right)  =\left(  x,0\right)  \text{; \ }Q\left(  x,y\right)
=\left(  0,y\right)  .
\]

Clearly, $P\circ Q=0$ and $P+Q=Id_{X\oplus_{N}Y}$.

The sum $X\oplus_{N}Y$ is said to be \textit{orthogonal} if $\left\|
P\right\|  =1$ (it will be constantly assumed that in an orthogonal sum
$X\oplus_{N}Y$ the norm one projection is associated with the first its
component). If $\left\|  P\right\|  =\left\|  Q\right\|  =1$ the sum
$X\oplus_{N}Y$ is said to be biorthogonal.

A particular case of a biorthogonal sum is an \textit{unconditional sum},
which is given by the condition%
\[
\left\|  P-Q\right\|  =1.
\]

The following result due to the folklore.

\begin{proposition}
Every unconditional sum $Z=X\oplus_{N}Y$ is biorthogonal.
\end{proposition}

\begin{proof}
Notice that $\left\|  P^{2}\right\|  \leq\left\|  P\right\|  \left\|
P\right\|  $; $\left\|  Q^{2}\right\|  \leq\left\|  Q\right\|  \left\|
Q\right\|  $ and, hence,
\[
\min\{\left\|  P\right\|  ,\left\|  Q\right\|  \}\geq1
\]
From the other hand,
\[
1=\left\|  P-Q\right\|  =\left\|  Id_{Z}-2Q\right\|  \geq\left\|  2Q\right\|
-\left\|  Id_{Z}\right\|  =2\left\|  Q\right\|  -1
\]
and, hence, $\left\|  Q\right\|  \leq1$. Similarly, $\left\|  P\right\|
\leq1$. So, $\left\|  P\right\|  =\left\|  Q\right\|  =1$.
\end{proof}

A particular case of an unconditional sums presents $l_{p}$\textit{-sum}
$X\oplus_{p}Y$, $p\geq1$, where the norm $N_{p}$ on $X\dotplus Y$ is given by
\begin{align*}
N_{p}\left(  x,y\right)   &  =\left(  \left\|  x\right\|  ^{p}+\left\|
y\right\|  ^{p}\right)  ^{1/p},\text{ }p\in\lbrack1,\infty)\text{; \ }\\
N_{\infty}\left(  x,y\right)   &  =\max\{\left\|  x\right\|  ,\left\|
y\right\|  \}.
\end{align*}

Notice that all norms $N\left(  \cdot,\cdot\right)  $ define on $X\dotplus Y$
the same topology: all spaces of kind $X\oplus_{N}Y$ are pairwice isomorphic.
This is not the case for infinite direct sums. The sum $\sum
\nolimits_{i<\infty}\oplus X_{i}$ (here the sign $\oplus$ is not related to
any concrete norm $N\left(  \cdot,\cdot\right)  $) is defined as the closure
of union of the chain $Z_{1}\hookrightarrow Z_{2}\hookrightarrow...$ where
$Z_{1}=X_{1}$; $Z_{n}=Z_{n-1}\oplus_{N_{n-1}}X_{n}$, all norms $\{N_{n}\left(
\cdot,\cdot\right)  :n<\infty\}$ are orthogonal and each embedding of
$Z_{n-1}$ into $Z_{n}$ is canonical.

The simplest example of infinite direct sums are $l_{p}$-sums (or shortly
$p$-sums); $1\leq p\leq\infty$, where\ all norms $N_{n}$ are just the $N_{p}$
norm. Such $p$-sum will be denoted by $\left(  \sum\nolimits_{i<\infty}\oplus
X_{i}\right)  _{p}$ or, simply, by $l_{p}\left(  X_{i}\right)  $. The space
$l_{p}\left(  X_{i}\right)  $ may be regarded as the space of sequences
$\frak{x}=\{x_{i}\in X_{i}:i<\infty\}$ with the finite norm
\[
\left\|  \frak{x}\right\|  _{p}=(\sum\nolimits_{p<\infty}\left\|
x_{i}\right\|  _{i}^{p})^{1/p}.
\]
When all spaces $\left(  X_{i}\right)  _{i<\infty}$ are equal to the given
Banach space $X$, their infinite $p$-sum is denoted by $l_{p}\left(  X\right)
$.

\begin{definition}
A class $X^{f}$ (and its Minkowski's base $\frak{M}(X^{f})$) is said to be
divisible if for every $A$, $B\in\frak{M}(X^{f})$ there exists an orthogonal
sum $C=A\oplus_{N}B$, which belongs to $\frak{M}(X^{f})$ as well.
\end{definition}

Certainly, every divisible class $X^{f}$ is generated by an
infinite-dimensional Banach space $X$. Indeed, for every $A\in\frak{M}(X^{f})$
spaces $A\oplus A$, $A\oplus A\oplus A$ etc. belong to $\frak{M}(X^{f})$ as
well (here and below when the lower index in $\oplus_{N}$ is omitted, the
symbol $\oplus$ denotes an arbitrary orthogonal sum).

The following result characterizes divisible classes.

\begin{theorem}
Let $X\in\mathcal{B}$; $X^{f}$ be the corresponding class of finite
equivalence. The following conditions are equivalent:

\begin{enumerate}
\item $X^{f}$ is divisible;

\item  There exists a space $W\in X^{f}$ such that an orthogonal sum $W\oplus
W$ also belongs to $X^{f}$;

\item  There exists $W\in X^{f}$ such that an orthogonal infinite sum
$\sum\oplus W_{i}$, where $W_{i}=W$, belongs to $X^{f}$

\item  For every $Z<_{f}X$ there exists such orthogonal sum that $Z\oplus Z$
is finitely representable in $X$ as well;

\item  For every $W\in X^{f}$ there exists an infinite orthogonal sum
$\sum\oplus W_{i}$, where $W_{i}=W$, which belongs to $X^{f}$;

\item  There exists $Y\in X^{f}$, which is finitely representable in each its
subspace of finite codimension;

\item  Every $Y\in X^{f}$ is finitely representable in each its subspace of
finite codimension.
\end{enumerate}
\end{theorem}

\begin{proof}
Surely, $\left(  4\Rightarrow2\right)  $; $\left(  5\Rightarrow3\right)  $ and
$\left(  7\Rightarrow6\right)  $.

Recall that every class $X^{f}$ is closed under direct limits of isometric
direct systems: if $\left\langle A_{i},m_{ij}\right\rangle _{i,j\in I}$ is a
direct system, all $m_{ij}:A_{i}\rightarrow A_{j}$ are isometric embeddings
and $A_{i}<_{f}X$ ($i\in I$) then $\underset{\rightarrow}{\lim}A_{i}$ also is
finitely representable in $X$.

$\left(  1\Rightarrow2\right)  $. Assume that $X^{f}$ is divisible. Let
$\left(  A_{i}\right)  $ be dense in $\frak{M}(X^{f}).$ Let $\oplus$ denotes
any direct sum that keep the class $\frak{M}(X^{f})$: if $A$, $B\in
\frak{M}(X^{f})$ then $A\oplus B$ is also assumed to be in the class
$\frak{M}(X^{f})$ too. Put
\begin{align*}
Z_{1}  &  =A_{1}\oplus A_{2};\text{ \ \ \ \ \ \ \ }Y_{1}=Z_{1}\oplus Z_{1};\\
Z_{n}  &  =Z_{n-1}\oplus A_{n+1};\text{ \ }Y_{n}=Z_{n}\oplus Z_{n}\text{
\ (}n>1\text{).}%
\end{align*}

Clearly, $W=\underset{\rightarrow}{\lim}Z_{n}\in X^{f}$ and $\underset
{\rightarrow}{\lim}Y_{n}=W\oplus W\in X^{f}$.

$\left(  2\Rightarrow3\right)  $. Proceeding by induction, it may be obtained
a sequence
\[
W_{1}=W;\text{ }W_{2}=W\oplus W;\text{ }...;\text{ }W_{n}=W_{n-1}\oplus
W;\text{ ...,}%
\]
which direct limit is a direct sum $\sum\oplus W\in X^{f}$.

$\left(  2\Rightarrow4\right)  $. Notice that ultrapowers preserve direct sums:

\textit{If }$Y=U\oplus V$\textit{ then there exists a such sum }%
$\oplus^{\prime}$ \textit{(that depends on the ultrafilter }$D$\textit{) that
}$\left(  Y\right)  _{D}=\left(  U\right)  _{D}\oplus^{\prime}\left(
V\right)  _{D}$\textit{.}

This observation easily follows from the standard properties of ultrapowers
(cf. [13]). Since every space that is finitely representable in $X$ is
isometric to a subspace of some ultrapower $\left(  X\right)  _{D}$, it is
clear that for every $Z<_{f}X$ it may be found such ultrafilter $D$ that
$\left(  W\right)  _{D}\oplus^{\prime}\left(  W\right)  _{D}$ contains a
subspace $Z\oplus^{\prime}Z$ (the space $W$ was defined above).

$\left(  3\Rightarrow5\right)  $. Let $Z<_{f}X$.

Then $Z\oplus^{\prime}Z\hookrightarrow\left(  W\right)  _{D}$ and
$Z\oplus^{\prime}Z\oplus^{\prime}Z\oplus^{\prime}Z\hookrightarrow\left(
W\right)  _{D}\oplus^{\prime}\left(  W\right)  _{D}$.

Proceed by induction, put $Z_{1=}Z\oplus^{\prime}Z$; $Z_{n+1}=Z_{n}%
\oplus^{\prime}Z_{n}$ for all $n\in\mathbb{N}$.

Surely, the direct limit of the chain $Z_{1}\hookrightarrow Z_{2}%
\hookrightarrow...\hookrightarrow Z_{n}\hookrightarrow...$is the desired space.

$(3\Rightarrow6)$. Clearly, any space of kind $\sum\oplus W$ (where the sum is
infinite) is finitely representable in each its subspace of finite codimension.

$\left(  6\Rightarrow1\right)  $. Let $Y\in X^{f}$ and $Y$ is finitely
representable in each its subspace of finite codimension. Clearly, $H(Y)$ is
dense in $\frak{M}(X^{f})$. Let $A$ be a finite dimensional subspace of $Y$.
According to [14], for every $\varepsilon>0$ there exists a subspace $B$ of
$Y$ of finite codimension, which is $(1+\varepsilon)$ - orthogonal to $A$,
i.e.\ such that $B\cap A=\varnothing$ and there exists a projections from
$\operatorname{span}\{A\cup B\}\hookrightarrow H(Y)$ onto $A$ of norm
$\left\|  P\right\|  \leq1+\varepsilon$. Hence, the property $5$ implies that
for every $A$, $C\in H(Y)$ (not necessary distinct) and for every
$\varepsilon>0$ some direct sum $A\oplus C\in H(Y)$ is $(1+\varepsilon)$ -
orthogonal, i.e., there exists a projection $P^{\prime}:A\oplus C\rightarrow
A$ with the norm $\left\|  P^{\prime}\right\|  \leq1+\varepsilon$. Since
$\varepsilon$ is arbitrary, and $H(Y)$ is dense in $\frak{M}(X^{f})$, it is
clear that this is equivalent to $1$.

$\left(  6\Rightarrow7\right)  $. Assume that for some $Y\in X^{f}$ and for
some its subspace $Y_{\left(  n\right)  }$ of $\operatorname{codim}Y_{\left(
n\right)  }=n<\infty$ there is a subspace $A\hookrightarrow Y$ such that for
some $\varepsilon>0$
\[
\inf\{d\left(  A,B\right)  :B\hookrightarrow Y_{\left(  n\right)  }%
\}\geq1+\varepsilon.
\]
This means that there exists such $n_{0}$ that for all $n>n_{0}$ the set
$\frak{M}(X^{f})$ does not contain any space of kind $A^{n}=A\oplus
A\oplus...\oplus A$ ($n$ times)$.$ This yields that $X^{f}$ is not divisible
and, hence, there no spaces of $X^{f}$ having the property $6$.
\end{proof}

The following result is obvious and its proof will be omitted.

\begin{theorem}
Any Banach space $X$ may be symmetrically embedded into a space
\[
l_{2}(X)=(\sum\nolimits_{i<\infty}\oplus X_{i})_{2},
\]
where all $X_{i}$ are isometric to $X$. Immediately, $l_{2}(X)$ generates a
divisible class $\mathsf{D}_{2}(X^{f})=\left(  l_{2}(X)\right)  ^{f}$ which is
of the same type and cotype as $X^{f}$ and is superreflexive if and only if
$X^{f}$ is.
\end{theorem}

\begin{remark}
The procedure $\mathsf{D}_{2}:X^{f}\rightarrow\left(  l_{2}(X)\right)  ^{f}$
may be regarded as a closure operator on the partially ordered set $f\left(
\mathcal{B}\right)  $ of all\ different classes of finite equivalence (which
may be considered to be a set of all different subsets of $\frak{M}$ of kind
$\frak{M}(X^{f})$, ordered by the standard relation $\subset$ of
set-theoretical inclusion)

Indeed, it is

\begin{itemize}
\item  Monotone, i.e., $X^{f}<_{f}\mathsf{D}_{2}(X^{f})$;

\item  Idempotent, i.e., $\mathsf{D}_{2}(X^{f})=\mathsf{D}_{2}(\mathsf{D}%
_{2}(X^{f}))$;

\item  Preserve the order: $X^{f}<_{f}Y^{f}$ $\Longrightarrow\mathsf{D}%
_{2}(X^{f})<_{f}\mathsf{D}_{2}(Y^{f})$.
\end{itemize}

It is of interest\ to notice that extreme points of $f\left(  \mathcal{B}%
\right)  $ - the class $\left(  l_{2}\right)  ^{f}$ (minimal with respect to
the relation $<_{f}$ by the Dvoretzki theorem) and the class $\left(
l_{\infty}\right)  ^{f}$ (maximal by an easy consequence of the Hahn-Banach
theorem) are stable under the procedure $\mathsf{D}_{2}$:
\[
\mathsf{D}_{2}(\left(  l_{2}\right)  ^{f})=\left(  l_{2}\right)  ^{f};\text{
}\mathsf{D}_{2}(\left(  c_{0}\right)  ^{f})=\left(  c_{0}\right)  ^{f}.
\]
\end{remark}

To distinguish between general divisible classes and classes of type
$\mathsf{D}_{2}(X^{f})$, the last ones will be called $\mathit{2}%
$\textit{-divisible classes}.

\section{Quotient-closed classes of finite equivalence}

\begin{definition}
A class $X^{f}$ of finite equivalence is said to be quotient-closed if for
every $A\in\frak{M}(X^{f})$ and each its subspace $B$ the quotient $A/B$
belongs to $\frak{M}(X^{f})$.
\end{definition}

As it will be shown, quotient-closed classes plays a significant role in the
theory of injective Banach spaces, defined in the next section.

In this section it will be shown how to enlarge a Minkowski's base
$\frak{M}(X^{f})$ of a certain $B$-convex class $X^{f}$ to obtain a set
$\frak{N}$, which forms a Minkowski's base $\frak{M}(W^{f})$ for some
quotient-closed class $W^{f}$. In other words, the main result of this section
is the proof of existence of wide collection of quotient-closed classes.

Let $K\subseteq\frak{M}$ be a class of finite-dimensional Banach spaces
(recall, that isometric spaces are identified). Define operations $H$, $Q$ and
$^{\ast}$ that transform a class $K$ to another class of finite-dimensional
Banach spaces - $H(K)$; $Q(K)$ or $(K)^{\ast}$ respectively. Namely, let%
\begin{align*}
H(K)  &  =\{A\in\frak{M}:A\hookrightarrow B;\text{ \ }B\in K\}\\
Q(K)  &  =\{A\in\frak{M}:A=B/F;\text{ \ }F\hookrightarrow B;\text{ \ }B\in
K\}\\
(K)^{\ast}  &  =\{A^{\ast}\in\frak{M}:A\in K\}
\end{align*}

In words, $H(K)$ consists of all subspaces of spaces from $K$; $Q(K)$ contains
all quotient spaces of spaces of $K$; $(K)^{\ast}$ contains all conjugates of
spaces of $K$.

The following theorem lists properties of these operations. In iteration of
the operations parentheses may be omitted.

Thus, $K^{\ast\ast}\overset{\operatorname*{def}}{=}\left(  (K)^{\ast}\right)
^{\ast}$; $HH(K)\overset{\operatorname*{def}}{=}H(H(K))$ and so on.

\begin{theorem}
Any set $K$ of finite-dimensional Banach spaces has the following properties:

\begin{enumerate}
\item $K^{\ast\ast}=K$; $HH(K)=H(K)$; $QQ(K)=Q(K)$;

\item $K\subset H(K)$; $K\subset Q(K)$;

\item  If $K_{1}\subset K_{2}$ then $H(K_{1})\subset H(K_{2})$ and
$Q(K_{1})\subset Q(K_{2})$;

\item $\left(  H(K)\right)  ^{\ast}=Q(K^{\ast})$; $\left(  Q(K)\right)
^{\ast}=H(K^{\ast})$;

\item $HQ(HQ(K))=HQ(K)$; $QH(QH(K))=QH(K)$.
\end{enumerate}
\end{theorem}

\begin{proof}
1, 2 and 3 are obvious.

4. If $A\in Q(K)$ then $A=B/E$ for some $B\in K$ and its subspace $E$. Thus,
$A^{\ast}$ is isometric to a subspace of $B^{\ast}.$ Hence, $A^{\ast}\in
H(B^{\ast})$, i.e., $A^{\ast}\in H(K^{\ast})$. Since $A$ is arbitrary,
$\left(  Q(K)\right)  ^{\ast}\subseteq H(K^{\ast})$. Analogously, if $B\in K$
and $A\in H(B)$ then $A^{\ast}$ may be identified with a quotient $B^{\ast
}/A^{\perp}$, where $A^{\perp}$ is the annihilator of $A$ in $B^{\ast}$:%
\[
A^{\perp}=\{f\in B^{\ast}:\left(  \forall a\in A\right)  \text{ \ }\left(
f\left(  a\right)  =0\right)  \}.
\]

Hence, $A^{\ast}\in(Q(K^{\ast}))^{\ast}$, and thus $\left(  H(K)\right)
^{\ast}\subseteq Q(K^{\ast})$.

From the other hand,%
\begin{align*}
H(K^{\ast}) &  =\left(  H(K^{\ast})\right)  ^{\ast\ast}\subseteq\left(
Q(K^{\ast\ast})\right)  ^{\ast}=\left(  Q(K)\right)  ^{\ast};\\
Q(K^{\ast}) &  =\left(  Q(K^{\ast})\right)  ^{\ast\ast}\subseteq\left(
H(K^{\ast\ast})\right)  ^{\ast}=\left(  H(K)\right)  ^{\ast}.
\end{align*}

5. Let $A\in HQ(K)$. Then $A$ is isometric to a subspace of some quotient
space $E/F$, where $E\in K$; $F\hookrightarrow E$. If $B$ is a subspace of $A$
then $\left(  A/B\right)  ^{\ast}=\left(  E/F\right)  ^{\ast}/B^{\perp}$, i.e.
$\left(  Q\left(  HQ(K)\right)  \right)  ^{\ast}\subseteq Q\left(  Q(K)^{\ast
}\right)  $. Because of
\begin{align*}
Q\left(  HQ(K)\right)   &  =\left(  Q\left(  HQ(K)\right)  \right)  ^{\ast
\ast}\subseteq(Q\left(  Q(K)^{\ast}\right)  )^{\ast}\\
&  \subseteq H(Q(K)^{\ast\ast})=HQ(K),
\end{align*}

we have%
\[
H(Q(H(Q(K))))\subseteq H(H(Q(K)))=HQ(K).
\]

Analogously, if $A\in QH(K)$ then $A$ is isometric to a quotient space $F/E$,
where $F\in H(B)$ for some $B\in K$ and $E\hookrightarrow B$. If $W\in H(A)$,
i.e., if $W\in H(F/E)$ then $W^{\ast}=(F/E)^{\ast}/W^{\perp}$ and
$(F/E)^{\ast}$ is isometric to a subspace $E^{\perp}$ of $F^{\ast}%
\in(H(B))^{\ast}$. Thus, $(H(QH(K)))^{\ast}\subseteq H((H(K))^{\ast})$ and
\begin{align*}
H\left(  QH(K)\right)   &  =\left(  H\left(  QH(K)\right)  \right)  ^{\ast
\ast}\subseteq(H\left(  H(K)^{\ast}\right)  )^{\ast}\\
&  \subseteq Q(H(K)^{\ast\ast})=QH(K).
\end{align*}

Hence,
\[
Q(H(Q(HQ(K))))\subseteq Q(Q(H(K)))=QH(K).
\]

Converse inclusion follows from 2.
\end{proof}

\begin{theorem}
Let $X$ generates the divisible class $X^{f}$. Consider Minkowski's base
$\frak{M}(X^{f})$ and its enlargement $H(Q(\frak{M}(X^{f})))=HQ\frak{M}%
(X^{f})$. Then there exists a Banach space $W$ such that $HQ\frak{M}%
(X^{f})=\frak{M}(W^{f})$.
\end{theorem}

\begin{proof}
Obviously, $HQ\frak{M}(X^{f})$ has properties (\textbf{H}) and (\textbf{C}).

Since $\frak{M}(X^{f})$ is divisible, then for any $A,B\in\frak{M}(X^{f})$ the
space $A\oplus B$ belongs to $\frak{M}(X^{f})$ and, hence, to $HQ\frak{M}%
(X^{f})$. If $A,B\in Q\frak{M}(X^{f})$ then $A=F/F_{1}$; $B=E/E_{1}$ for some
$E,F\in Q\frak{M}(X^{f})$.

$F/F_{1}\oplus E$ is isometric to a space ($F\oplus E)/F_{1}^{\prime}$, where
\[
F_{1}^{\prime}=\{(f,0)\in F\oplus E:f\in F_{1}\}
\]
and, hence, belongs to $Q\frak{M}(X^{f})$. Thus,
\[
F/F_{1}\oplus E/E_{1}=(F/F_{1}\oplus E)/E_{1}^{\prime},
\]
where
\[
E_{1}^{\prime}=\{(o,e)\in F/F_{1}\oplus E:e\in E_{1}\}
\]
and, hence, belongs to $Q\frak{M}(X^{f})$ as well.

If $A,B\in HQ\frak{M}(X^{f})$ then $A\hookrightarrow E$, $B\hookrightarrow F$
for some $E,F\in Q\frak{M}(X^{f})$.

$E\oplus F\in Q\frak{M}(X^{f})$ and, hence, $A\oplus B\in HQ\frak{M}(X^{f})$.

Thus, $HQ\frak{M}(X^{f})$ has the property (\textbf{A}$_{0}$). The desired
result follows from the preceding theorem.
\end{proof}

Define on $f(\mathcal{B})$ the procedure $\ast:f\left(  \mathcal{B}\right)
\rightarrow f\left(  \mathcal{B}\right)  $. Let $Y\in\mathcal{B}$; $Y^{f}$ be
a divisible class; $\frak{M}(Y^{f})$ - its Minkowski's base. Clearly,
$H\frak{M}(Y^{f})=\frak{M}(Y^{f})$. Hence,
\[
(\frak{M}(Y^{f}))^{\ast}=H(\frak{M}(Y^{f}))^{\ast}=Q(\frak{M}(Y^{f}))^{\ast}%
\]
and we can put $\frak{M}(\ast Y^{f})\overset{\operatorname{def}}{=}%
H(\frak{M}(Y^{f}))^{\ast}$. From the preceding formulae it follows that
$\frak{M}(\ast Y^{f})=HQ\frak{M}(\ast Y^{f})$.$\ $

From divisibility of $Y^{f}$ it follows that $H(\frak{M}(Y^{f}))^{\ast
}=HQ(\frak{M}(Y^{f}))^{\ast}$ satisfies the condition $\left(
\text{\textbf{A}}_{0}\right)  $. Hence $H(\frak{M}(Y^{f}))^{\ast}$ is the
Minkowski's base of some class of finite equivalence, which will be denoted by
$\ast Y^{f}$. Clearly, the class $\ast Y^{f}$ is quotient-closed and divisible.

In a case, when $Y^{f}$ is not divisible it will be considered the class
$X^{f}=\mathsf{D}_{2}Y^{f}$.

\begin{definition}
The procedure $\ast:f\left(  \mathcal{B}\right)  \rightarrow f\left(
\mathcal{B}\right)  $ is given by:%
\begin{align*}
&  \text{ \ \ \ \ If }Y^{f}\text{ is divisible then }\frak{M}(\ast
Y^{f})\overset{\operatorname{def}}{=}H(\frak{M}(Y^{f}))^{\ast};\\
&  \text{If }X^{f}\text{ is not divisible then }\frak{M}(\ast X^{f}%
)\overset{\operatorname{def}}{=}H(\frak{M}(\mathsf{D}_{2}X^{f}))^{\ast}.
\end{align*}
\end{definition}

The procedure $\ast$ may be also defined in other way.

Let $Y\in\mathcal{B}$; $Y^{f}$ be a divisible class. Let $\left(
Y_{n}\right)  _{n<\infty}$ be a countable dense subset of $\frak{M}(Y^{f})$.
Consider a space $Z_{\oplus}=\sum_{n<\infty}\oplus Y_{n}$ and its conjugate
$Z_{\oplus}^{\ast}$.

$Z_{\oplus}^{\ast}$ generates a class $\left(  Z_{\oplus}^{\ast}\right)  ^{f}$
which will be regarded as a result of acting of the procedure $\ast
:Y^{f}\rightarrow\left(  Z_{\oplus}^{\ast}\right)  ^{f}$. Iterations of the
procedure $\ast$ are given by following steps.

Let $\left(  Z_{n}\right)  _{n<\infty}$ be a countable dense subset of
$\frak{M}(\left(  Z_{\oplus}^{\ast}\right)  ^{f})$. Consider a space
$W=\sum_{n<\infty}\oplus Z_{n}$ and its conjugate $W^{\ast}$.
Clearly,\ $W^{\ast}$ generates a class
\[
\left(  W^{\ast}\right)  ^{f}=\ast\left(  Z_{\oplus}^{\ast}\right)  ^{f}%
=\ast\ast(Y^{f}).
\]

Another definition of the 'double star' procedure may be given by an
enlargement of the Minkowski's base.

\begin{definition}
The procedure $\ast\ast:f\left(  \mathcal{B}\right)  \rightarrow f\left(
\mathcal{B}\right)  $ is given by:%
\begin{align*}
&  \text{ \ \ \ \ If }Y^{f}\text{ is divisible then }\frak{M}(\ast\ast
Y^{f})\overset{\operatorname{def}}{=}HQ(\frak{M}(Y^{f}));\\
&  \text{If }X^{f}\text{ is not divisible then }\frak{M}(\ast\ast
X^{f})\overset{\operatorname{def}}{=}HQ(\frak{M}(\mathsf{D}_{2}X^{f}))
\end{align*}
\end{definition}

Certainly, for every $X\in\mathcal{B}$, $X^{f}<_{f}\ast\ast X^{f}$. Let us
show that the procedure $\ast\ast$ holds the type, cotype and superreflexivity.

\begin{theorem}
Let $X$ be a Banach space, which generates a class of finite equivalence
$X^{f}$. If $S(X^{f})$ is its $l_{p}$-spectrum then the procedure $\ast\ast$
maps $X^{f}$ to a class $\ast\ast(X^{f})$ of the upper and lower bounds of
$S(\ast\ast(X^{f}))$ as $X^{f}$. In difference from $S(X^{f})$ it is always
the closed segment
\[
S(\ast\ast(X^{f}))=\left[  \inf S(X^{f});\sup S(X^{f})\right]  .
\]
If $X$ is superreflexive then $\ast\ast(X^{f})$ is superreflexive too.
\end{theorem}

\begin{proof}
From the construction of the class $\ast\ast(X^{f})$ it follows that for every
$p\in S(X^{f})$ the number $q^{\prime}=p/(p-1)\in S(\ast(X^{f}))$ and
conversely, $q\in S(\ast(X^{f}))$ implies that $p^{\prime}=q/(q-1)\in
S(X^{f})$. It is clear that the procedure $\ast$ holds superreflexivity. From
the preceding theorem it follows that $\ast\ast$ enjoys this property as well.
The invariance of type and cotype follows from results of G. Pisier [15]
\end{proof}

\begin{remark}
If $X$ is not $B$-convex, then
\[
\ast X^{f}=\ast\ast X^{f}=\left(  c_{0}\right)  ^{f}.
\]
\end{remark}

\begin{theorem}
For any Banach space $X$ the class $\ast\ast X^{f}$ is divisible.
\end{theorem}

\begin{proof}
Let $\frak{N}=\frak{M}(\ast\ast(X^{f}))$.

Since for any pair $A,B\in\frak{N}$ \ their direct orthogonal sum belongs to
$\frak{N}$, by induction, $\sum\nolimits_{i\in I}\oplus A_{i}\in\frak{N}$ for
any finite subset $\{A_{i}:i\in I\}\subset\frak{N}$.

Hence, any infinite direct sum $\sum\nolimits_{i\in I}\oplus A_{i}$ is
finitely representable in $\ast\ast X^{f}$.

Let $\{A_{i}:i<\infty\}\subset\frak{N}$ be dense in $\frak{N}$.

Let $Y_{1}=\sum\nolimits_{i<\infty}\oplus A_{i}$; $Y_{n+1}=Y_{n}\oplus Y_{1}$;
$Y_{\infty}=\overline{\cup Y_{n}}$ (the upper line denotes the closure).

Clearly, $Y_{\infty}$ belongs to $\ast\ast X^{f}$.
\end{proof}
\end{document}